\documentclass[12pt,twocolumn]{amsart}
\usepackage{a4}

\usepackage[latin1]{inputenc}
\usepackage{enumerate}
\usepackage{hyperref}
\usepackage{graphicx}
\usepackage{amssymb}
\usepackage{xypic}
\usepackage[vmargin=2cm,hmargin=2cm]{geometry}
\usepackage{color}

\makeindex

\def\N{{\mathbb N}}
\def\Q{{\mathbb Q}}
\def\Z{{\mathbb Z}}

\def\Ap[#1,#2]{{{\rm Ap}\left(#1,#2\right)}} %Ap\'ery sets [S, elements]

\def\setc[#1,#2]{\left\{  #1 ~|~ #2 \right\}} %cto \{#|#\}
\def\m{{\rm m}}
\def\e{{\rm e}}
\def\F{{\rm F}}
\def\H{{\rm H}}
\def\g{{\rm g}}
\def\t{{\rm t}}
\def\PF{{\rm PF}}

\def\Over{{\mathcal O}}
\def\Sem{{\mathcal S}}

\def\D{{\rm D}}
\def\FH{{\rm FG}}

\def\S{{\rm S}}
\def\T{{\rm T}}
\def\P{{\rm P}}
\def\Si[#1]{{\S\left(\left[#1\right]\right)}}
\def\So[#1,#2]{{\S\left(\left]\frac{#2}{#1},\frac{#2}{#1-1}\right[\right)}}
\def\tb[#1,#2,#3]{{\frac{\langle #1, \, #2\rangle}{#3}}}
\def\Gr{{\mathcal{G}}}

\def\FV{{\mathcal F}}
\def\VV{{\mathcal V}}

\def\Free[#1]{{{\rm Free}(#1)}}
\def\Cong[#1]{{{\rm Cong}(#1)}}

\hypersetup{%
  pdftitle={Numerical semigroups problem list},%
  pdfauthor={M. Delgado, Pedro A. García Sánchez, J. C. Rosales}
}

\title{Numerical semigroups problem list}

\author{M. Delgado}
\email{mdelgado@fc.up.pt}

\address{CMUP and DM-FCUP, Universidade do Porto, Portugal}

\author{P. A. García-Sánchez}
\email{pedro@ugr.es}

\author{J. C. Rosales}
\email{jrosales@ugr.es}

\address{Departamento de Álgebra, Universidad de Granada, España}

\begin{document}

\twocolumn[
\begin{@twocolumnfalse}
\maketitle
\end{@twocolumnfalse}
]
\section{Notable elements and first problems}

A numerical semigroup is a subset of $\N$ (here $\N$ denotes the set of nonnegative integers) that is closed under addition, contains the zero element, and its complement in $\N$ is finite.

If $A$ is a nonempty subset of $\N$, we denote by $\langle A\rangle$ the submonoid of $\N$ generated by $A$, that is,
\[\langle A\rangle =\{\lambda_1 a_1+\cdots+ \lambda_n a_n~|~ n\in \N, \lambda_i\in \N, a_i\in A \}.\]
It is well known (see for instance \cite{libro,libro-sn}) that  $\langle A\rangle$ is a numerical semigroup if and only if $\gcd(A)=1$.

If $S$ is a numerical semigroup and $S=\langle A\rangle$ for some $A\subseteq S$, then we say that $A$ is a system of generators of $S$, or that $A$ generates $S$. Moreover, $A$ is a minimal system of generators of $S$ if no proper subset of $A$ generates $S$. In  \cite{libro-sn} it is shown that every numerical semigroup admits a unique minimal system of generators, and it has finitely many elements.

Let $S$ be a numerical semigroup and let $\{n_1<n_2<\cdots <n_e\}$ be its minimal system of generators. The integers $n_1$ and $e$ are known as the multiplicity and embedding dimension of $S$, and we will refer to them by using $\m(S)$ and $\e(S)$, respectively. This notation might seem amazing, but it is not so if one takes into account that there exists a large list of manuscripts devoted to the study of analytically irreducible one-dimensional local domains via their value semigroups, which are numerical semigroups. The invariants we just introduced, together with others that will show up later in this work, have an interpretation in that context, and this is why they have been named in this way. Along this line, \cite{barucci} is a good reference for the translation for the terminology used in the Theory of Numerical Semigroups and Algebraic Geometry.

Frobenius (1849-1917) during his lectures proposed the problem of giving a formula for the greatest integer that is not representable as a linear
  combination, with nonnegative integer coefficients, of a fixed set of integers with greatest common divisor equal to 1. He also raised the
question of determining how many positive integers do not admit such a representation. With our terminology, the first problem is equivalent to that
of finding a formula in terms of  the generators of a numerical semigroup $S$ of the greatest integer not belonging to $S$ (recall that its complement in $\N$ is finite). This number is thus known in the literature as the Frobenius number of $S$, and we will denote it by $\F(S)$. The elements of
     $\H(S)=\N\setminus S$ are called gaps of $S$. Therefore the second problem consists in determining the cardinality of $\H(S)$, sometimes known
as genus of $S$ (\cite{komeda}) or degree of singularity of $S$ (\cite{barucci}).

 In \cite{sylvester} Sylvester solves the just quoted  problems of Frobenius for embedding dimension two. For semigroups with embedding dimension greater than or equal to three these problems remain open. The current state of the problem is quite well collected in \cite{alfonsin}.

 Let $S$ be a numerical semigroup. Following the terminology introduced in  \cite{interseccion-sim} an integer $x$ is said to be a pseudo-Frobenius number  of $S$ if $x\not\in S$ and $x+S\setminus\{0\}\subseteq S$. We will denote by $\PF(S)$ the set of  pseudo-Frobenius numbers of $S$. The cardinality of $\PF(S)$ is called the type of $S$ (see \cite{barucci}) and we will denote it by $\t(S)$. It is proved in \cite{fgh} that if $\e(S)=2$, then $\t(S)=1$, and if $\e(S)=3$, then $\t(S)\in\{1,2\}$. It is also shown that if $\e(S)\geq 4$, then $\t(S)$ can be arbitrarily large, $\t(S)\leq \m(S)-1$ and that $(\t(S)+1)\g(S)\leq \t(s)(\F(S)+1)$. This is the starting point of a new line of research that consists in trying to determine the type of a numerical semigroup, once other invariants like multiplicity, embedding dimension, genus or Frobenius number are fixed.

Wilf in \cite{wilf} conjectures that if  $S$ is a numerical semigroup, then $\e(S)\g(S)\leq (\e(S)-1)(\F(S)+1)$. Some families of numerical semigroups for which it is known that the conjecture is true are collected in \cite{dobbs}. Other such families can be seen in \cite{kaplan,alessio}. The general case remains open.

Bras-Amorós computes in  \cite{b-a-conj} the number of numerical semigroups with genus $g\in\{0,\ldots,50\}$, and conjectures that the growth is similar to that of Fibonacci's sequence. However it has not been proved yet that there are more semigroups of genus  $g$ than of genus $g+1$. Several attempts already appear in the literature. Kaplan~\cite{kaplan} uses an approach that involves counting the semigroups by genus and multiplicity. He poses many related conjectures which could be taken literally and be posed here as problems. We suggest them to the reader.
A different approach, dealing with the asymptoptical behavior of the sequence of the number of numerical semigroups by genus, has been followed by Zhao~\cite{zhao}. Some progress has been achieved by Zhai~\cite{zhai}, but many questions remain open.  

\section{Proportionally modular semigroups}

 Following the terminology introduced in \cite{proportionally}, a proportionally modular Diophantine inequality   is an expression of the form $ax \bmod b\leq cx$, with $a$, $b$ and $c$ positive integers. The integers $a$, $b$ and $c$ are called the factor, the modulus and the proportion of the inequality, respectively. The set $\S(a,b,c)$ of solutions of the above inequality  is a numerical semigroup. We say that a numerical semigroup is proportionally modular if it is the set of solutions of some proportionally modular Diophantine inequality.

 Given a nonempty subset $A$ of $\Q_0^+$, we denote by $\langle A\rangle$ the submonoid of $(\Q_0^+,+)$ generated by $A$, whose definition is the same of that used in the previous section. Clearly, $\S(A)=\langle A\rangle \cap \N$ is a submonoid of $\N$. It is proved in \cite{proportionally} that if $a$, $b$ and $c$ are positive integers with $c<a<b$, then $\S(a,b,c)=\S\left(\left[ \frac{b}a, \frac{b}{a-c}\right]\right)$. Since $\S(a,b,c)=\N$ when $a\geq c$, and the inequality $ax\bmod b\leq cx$ has the same integer solutions as $(a\bmod b)x\bmod b\leq cx$, the condition
 $c<a<b$ is not restrictive.

 As a consequence of the results proved in \cite{proportionally}, we have that a numerical semigroup $S$ is proportionally modular if and only if there exist two positive rational numbers  $\alpha<\beta$ such that $S=\S([\alpha,\beta])$. This is also equivalent to the existence of an
 interval $I$, with nonempty interior, of the form $S=\S(I)$ (see \cite{opened}).

 By using the notation introduced in \cite{tres-anos}, a sequence of fractions $\frac{a_1}{b_1}<\frac{a_2}{b_2}<\cdots <\frac{a_p}{b_p}$ is said to be a Bézout sequence if $a_1,\ldots,a_p$, $b_1,\ldots, b_p$ are positive integers and $a_{i+1}b_i-a_ib_{i+1}=1$ for all $i\in\{1,\ldots, p-1\}$. The importance of the Bézout sequences in the study of proportionally modular semigroups highlights in the following result proved in \cite{tres-anos}. If $\frac{a_1}{b_1}<\frac{a_2}{b_2}<\cdots <\frac{a_p}{b_p}$ is a Bézout sequence, then $\S\left(\left[ \frac{a_1}{b_1},\frac{a_p}{b_p}\right]\right)=\langle a_1,\ldots, a_p\rangle$.

 A Bézout sequence $\frac{a_1}{b_1}<\frac{a_2}{b_2}<\cdots <\frac{a_p}{b_p}$ is proper if $a_{i+h}b_i-a_ib_{i+h}\geq 2$ for all $h\geq 2$ with $i,i+h\in \{1,\ldots,p\}$. Clearly, every Bézout sequence can be reduced (by removing some terms) to a proper Bézout sequence with the same ends as the original one. It is showed in  \cite{bullejos}, that if $\frac{a_1}{b_1}<\frac{a_2}{b_2}$ are two reduced fractions, then there exists an unique proper Bézout sequence with ends $\frac{a_1}{b_1}$ and $\frac{a_2}{b_2}$. Furthermore, in this work a procedure for obtaining this sequence is given.

 It is proved in \cite{tres-anos} that if $\frac{a_1}{b_1}<\frac{a_2}{b_2}<\cdots <\frac{a_p}{b_p}$ is a proper Bézout sequence, then there exists $h\in\{1,\ldots,p\}$ such that $a_1\geq \cdots \geq a_h\leq \cdots \leq a_p$ (the sequence $a_1,\ldots,a_p$ is convex). The following characterization is also proved there: a numerical semigroup is proportionally modular if and only if there exists a convex ordering if its minimal  generators $n_1,\ldots,n_e$ such that $\gcd\{n_i,n_{i+1}\}=1$ for all $i\in\{1,\ldots,e-1\}$ and $n_{j-1}+n_{j+1}\equiv 0\pmod {n_j}$ for all $j\in \{2,\ldots,e-1\}$.

A modular Diophantine inequality is a proportionally modular Diophantine inequality with proportion equal to one. A numerical semigroup is said to be modular if it is the set of solutions of some modular Diophantine inequality. Clearly, every modular numerical semigroup is proportionally modular, and this inclusion is strict as it is proved in \cite{proportionally}. A formula for $\g(\S(a,b,1))$ in function of $a$ and $b$ is given in \cite{modular}. The problems of finding formulas for $\F(\S(a,b,1))$, $\m(\S(a,b,1))$, $\t(\S(a,b,1))$ and $\e(\S(a,b,1))$ remain open. It is not known if the mentioned conjecture of Wilf is true for modular semigroups neither.

 A semigroup of the form $\{0,m,\to\}$ is said to be ordinary. A numerical semigroup $S$ is  an open modular numerical semigroup if it is ordinary or of it is  the form $S=\S\left(\left]\frac{b}a,\frac{b}{a-1}\right[\right)$ for some integers $2\leq a< b$. Therefore these semigroups are proportionally modular. Moreover, it is proved in \cite{opened} that every proportionally modular numerical semigroup can be expressed as a finite intersection of open modular numerical semigroups. The formulas $\F\left(\S\left(\left]\frac{b}a,\frac{b}{a-1}\right[\right)\right)$ and $\g\left(\S\left(\left]\frac{b}a,\frac{b}{a-1}\right[\right)\right)$ are also obtained in the just quoted work. However the rest of the problems previously suggested for modular numerical semigroups remain still open.

As we mentioned above, a characterization for proportionally modular numerical semigroups in terms of its systems of minimal generators is given in  \cite{tres-anos}. The question of giving formulas for the Frobenius number, genus and type of a proportionally modular numerical semigroup in terms of its system of minimal generators remains unsolved too.

 Following the terminology in \cite{contraidas}, a contracted modular Diophantine inequality is an expression of the form   $a x\bmod b\leq x-c$, where $a$, $b$ and $c$ are nonnegative integers and   $b\neq 0$.  Let us denote by $\T(a,b,c)$ the set of integer solutions of the last inequality. Then $\T(a,b,c)\cup\{0\}$ is a numerical semigroup. An algorithm that allows us to determine whether a semigroup is the set of solutions of a contracted modular Diophantine equation  is given in \cite{contraidas}. A formula for the genus of  $\T(a,b,c)\cup\{0\}$ is also given there.

A contracted proportionally modular Diophantine inequality is an expression of the form  $a x\bmod b \leq c x -d$, with $a$, $b$, $c$ and $d$ nonnegative integers and $b\neq 0\neq c$. If we denote by  $\P(a,b,c,d)$ the set of solutions of such an inequality, then $\P(a,b,c,d)\cup \{0\}$ is a numerical semigroup. It is not yet known an algorithm to determine whether a semigroup is of this form.

 The Stern-Brocot tree gives a recursive method for constructing all the reduced fractions $\frac{x}y$, with $x$ and $y$ positive integers (see \cite{gkp}). For constructing this tree we start with the expressions $\frac{0}1$ and $\frac{1}0$. In each step of the process we insert between each two consecutive expressions $\frac{m}n$ and $\frac{m'}{n'}$ its median $\frac{m+m'}{n+n'}$. We obtain in this way the sequences
 \[
\begin{array}{c}
  \frac{0}1<\frac{1}1<\frac{1}0\\ \\
  \frac{0}1<\frac{1}2<\frac{1}1<\frac{2}1<\frac{1}0\\ \\
  \frac{0}1<\frac{1}3<\frac{1}2<\frac{2}3<\frac{1}1<\frac{3}2<\frac{2}1<
  \frac{3}1<\frac{1}0\\
  \ldots
\end{array}
\]
 The Stern-Brocot tree can now be obtained by connecting each median with the fractions used for computing it and being in the previous level but not in the  levels above it.
{\tiny 
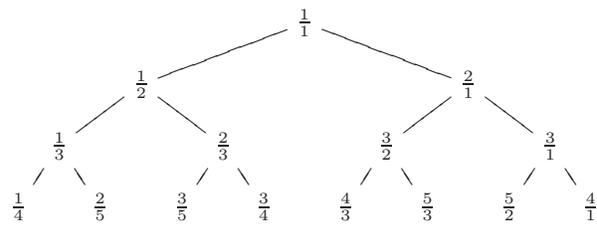
\xymatrix @R=0.5pc @C=.2pc{
 & & & & & & & {\frac{1}{1}}\ar@{-}[dllll]\ar@{-}[drrrr]& & & & & & & \\
 & & & {\frac{1}{2}}\ar@{-}[dll]\ar@{-}[drr] & & & & & & & &
 {\frac{2}{1}}\ar@{-}[dll]\ar@{-}[drr] & & & \\
 & {\frac{1}{3}}\ar@{-}[dl]\ar@{-}[dr] & & & &
 {\frac{2}{3}} \ar@{-}[dl]\ar@{-}[dr] & & & & {\frac{3}{2}}\ar@{-}[dl]\ar@{-}[dr] & & & & {\frac{3}{1}}\ar@{-}[dl]\ar@{-}[dr] &
 \\ \frac{1}{4} & & \frac{2}{5} & & \frac{3}{5} & & \frac{3}{4} & & \frac{4}{3} & & \frac{5}{3} & & \frac{5}{2} & & \frac{4}{1} }
}

 It is proved in \cite{bullejos} that if $\frac{m}n$ is the common  predecessor of two fractions $\frac{a}b< \frac{c}d$ in the Stern-Brocot tree, then $m$ is the multiplicity of $\S\left(\left[\frac{a}b,\frac{c}d\right]\right)$. It could be nice to obtain other constants of the semigroup by looking at this tree.

\section[Quotients]{The quotient of a numerical semigroup by a positive integer}

Let $S$ be a numerical semigroup and $p$ be a positive integer. Let us denote by
\[\frac{S}p=\{ x\in \N ~|~ px\in S\}.\]
Clearly, $\frac{S}p$ is a numerical semigroup, and we will call it the quotient of $S$ by $p$. According to this notation, we will call $\frac{S}2$ one half of $S$ and that $\frac{S}4$ is a quarter of $S$. These two cases will have an special importance in this section.

 It is proved in \cite{full} that a numerical semigroup  is proportionally modular if and only if it is the  quotient of an embedding dimension two numerical semigroup by a positive integer. This result is improved in \cite{aureliano} by proving  that a numerical semigroup  is proportionally modular if and only if it is of the form $\frac{\langle a,a+1\rangle}d$ with $a$ and $d$ positive integers. We still do not have formulas for $\F\left(\frac{\langle a,a+1\rangle}d\right)$, $\g\left(\frac{\langle a,a+1\rangle}d\right)$, $\m\left(\frac{\langle a,a+1\rangle}d\right)$, $\t\left(\frac{\langle a,a+1\rangle}d\right)$ and $\e\left(\frac{\langle a,a+1\rangle}d\right)$.

The next step in this line of research would be studying those numerical semigroups that are the quotient of a numerical semigroup with embedding dimension three by a positive integer. Unfortunately we do not have a procedure that allows us to distinguish such a semigroup from the rest. Moreover, we still do not know of any example of semigroups that are not of this form.

 A numerical semigroup $S$ is symmetric if $x\in \Z\setminus S$ implies $\F(S)-x\in S$. These semigroups have been widely studied. Their main motivation comes from a work by Kunz (\cite{kunz}) from which it can be deduced that  a numerical semigroup is symmetric if and only if its associated numerical semigroup ring is Gorenstein. Symmetric numerical semigroups always have odd Frobenius number, thus for numerical semigroups with even Frobenius number, the equivalent notion to symmetric semigroups is that of pseudo-symmetric numerical semigroups.  We say that $S$ is a pseudo-symmetric numerical semigroup if it has even Frobenius number and for all $x\in \Z\setminus S$, we have either $\F(S)-x\in S$ or $x=\frac{\F(S)}2$. The concept of irreducible semigroup, introduced in \cite{irreducibles-pacific},  collects these two families of semigroups. A numerical semigroup is irreducible if it cannot be expressed as the intersection of two semigroups that contain it properly. It can be
  proved that a semigroup is irreducible if and only if it is either symmetric (with odd Frobenius number) or pseudo-symmetric (with even Frobenius number).

 Intuition  (and the tables of the number of numerical semigroups with a given genus or Frobenius number we have) tells us that the percentage of irreducible numerical semigroups is quite small. It is proved in \cite{quotient-sym} that every numerical semigroup is one half of an infinite number of symmetric numerical semigroups. The apparent parallelism between symmetric and pseudo-symmetric numerical semigroups fails as we can see in \cite{quotient-psim}, where it is proved that a numerical semigroup is irreducible if and only it is one half of a pseudo-symmetric numerical semigroup. As a consequence we have that every numerical semigroup  is a quarter of  infinitely many pseudo-symmetric numerical semigroups. In \cite{swanson}, it is also shown that for every positive integer $d$ and every numerical semigroup $S$, there exist infinitely many symmetric numerical semigroups $T$ such that $S=T/d$, and if $d\ge 3$, then there exist infinitely many pseudo-symmetric numerical semigroups $T$ with $S=T/d$.

From the definition,  we deduce that a numerical semigroup $S$ is symmetric if and only if $\g(S)=\frac{\F(S)+1}2$. Therefore these numerical semigroups verify Wilf's conjecture previously mentioned. We raise the following question. If a numerical semigroup verifies Wilf's conjecture, then does so its half?

It can easily be seen that every numerical semigroup can be expressed as a finite intersection of irreducible numerical semigroups. A procedure for obtaining such a decomposition is given in \cite{oversemigroups}. Furthermore it is also explained how to obtain a decomposition with the least possible number of irreducibles. We still do not know how many numerical semigroups appear in these minimal decompositions, moreover, we wonder if there exists a positive integer $N$ such that every numerical semigroup can be expressed as an intersection of at most $N$ irreducible numerical semigroups.

 In \cite{toms} Toms introduces a class of numerical semigroups that are the positive cones of the $K_0$ groups of certain $C^*$-algebras. Given a numerical semigroup we say, inspired in this work, that it admits a Toms decomposition if and only if there exist positive integers  $q_1,\ldots, q_n$, $m_1,\ldots,m_n$ and $L$ such that $\gcd\{q_i,m_i\}= \gcd \{L,m_i\}= \gcd\{L,q_i\}=1$ for all $i\in\{1,\ldots,n\}$ and $S=\frac{1}L\bigcap_{i=1}^n\langle q_i,m_i\rangle$.

As $\frac{1}L\bigcap_{i=1}^n\langle q_i,m_i\rangle = \bigcap_{i=1}^n\frac{\langle q_i,m_i\rangle}L$, we have that  if a numerical semigroup admits a Toms decomposition, then  $S$ is a finite intersection of proportionally modular numerical semigroups. It is proved in \cite{ns-toms} that the reciprocal is also true. Therefore, a numerical semigroup admits a Toms decomposition if and only if it is an intersection of finitely many proportionally modular numerical semigroups. These kind of semigroups are studied in  \cite{spm}, where  an algorithm for distinguishing whether a numerical semigroup is an intersection of finitely many proportionally modular numerical semigroups  is given. Furthermore, in the affirmative case it gives us a minimal decomposition, and in the negative case it gives us the least numerical semigroup which is intersection of proportionally modular semigroups and contains the original numerical semigroup (its proportionally modular closure).

It is conjectured in \cite{contraidas} that every contracted modular numerical semigroup admits a Toms decomposition.

 Note that the numerical semigroups that admit a Toms decomposition are those that are the set of solutions of a system of proportionally modular Diophantine inequalities. It is proved in \cite{aureliano} that two systems of inequalities are always equivalent to another  system with all
 the inequalities having the same modulus, which moreover can be chosen to be prime. Now we raise the following question: is every system of proportionally modular Diophantine inequalities equivalent to a system with all proportions being equal to one?, or equivalently, if a numerical
 semigroup admits a Tom decomposition, can it be expressed as an intersection of  modular numerical semigroups?

Following the terminology introduced in  \cite{fg}, a gap $x$ in a numerical semigroup  $S$  is said to be fundamental if  $\{2x,3x\}\subset S$ (and therefore  $kx\in S$ for every integer with $k\geq 2$). Let us denote by $\FH(S)$ the set of all fundamental gaps of  $S$. If $X\subseteq \Z$, then $\D(X)$ will denote  the union of all positive divisors of the elements of $X$. It can easily be shown that  $S=\N\setminus \D(\FH(S))$. Therefore, a way to represent a semigroup is by giving its fundamental gaps. This representation is specially useful when studying the quotient of a semigroup  $S$ by a positive integer $d$, since $\FH\left(\frac{S}d\right)=\left\{ \frac{h}d ~|~h\in \FH(S),\ h\equiv 0\pmod d \right\}$.

The cardinality of the set of fundamental gaps of a semigroup is an invariant of the semigroup. We can therefore open a new line of research by studying numerical semigroups attending to their number of fundamental gaps. It would be also interesting to find simple sufficient conditions  that allow us to decide when a subset  $X$ of $\N$ is the set of fundamental gaps of some numerical semigroup.

 Let $S$ be a numerical semigroup. In \cite{dobles} the set $T$ of all numerical semigroups such that    $S=\frac{T}2$ is studied, the semigroup of the ``doubles'' of $S$. In the just quoted work we raise the question of finding a formula that depends on   $S$ and allows us to compute the minimum of the Frobenius numbers of the doubles of  $S$.

 Following this line we can ask ourselves about the set of all  ``triples'' (or multiples in  general) of a numerical semigroup.

Finally, it would be interesting to characterize the families of numerical semigroups verifying that any of its elements can be realized as a quotient of some element of the family by a fixed positive integer.

%\textcolor{blue}{Trabajo de Dobbs y su alumno \\ Trabajo de Swanson}

\section{Frobenius Varieties}

A directed graph $G$ is a pair $(V,E)$, where $V$ is a nonempty set whose elements are called vertices, and $E$ is a subset of  $\{(u,v)\in V\times V~|~ u\neq v\}$. The elements of $E$ are called edges of the graph. A path connecting two vertices  $x$ and $y$ of $G$ is a sequence of distinct edges of the form $(v_0,v_1), (v_1,v_2), \ldots, (v_{n-1},v_n)$ with $v_0=x$ and $v_n=y$. A graph $G$ is a tree if there exists a vertex  $r$ (called the root of $G$) such that for any other vertex  $x$ of $G$, there exists an unique path connecting $x$ and $r$. If $(x,y)$ is an edge of the tree, then $x$ is a son of  $y$. A vertex of a tree is a leaf if it has no sons.

 Let $\Sem$ be the set of all numerical semigroups. We define the graph associated to $\Sem$, $\Gr(\Sem)$, to be the graph whose vertices are all the elements of  $\Sem$ and $(T,S)\in \Sem\times \Sem$ is an edge if $S=T\cup \{\F(T)\}$. In \cite{libro-sn}, it is proved that $\Gr(\Sem)$ is a tree with root $\N$, and that the sons of $S\in \Sem$ are the subsets $S\setminus\{x_1\},\ldots, S\setminus\{x_r\}$, where $x_1,\ldots,x_r$ are the minimal generators of $S$ greater than $\F(S)$. Therefore  $S$ is a leaf of $\Gr(\Sem)$ if it has no minimal generators greater than $\F(S)$. These results allow us to construct recursively the set of numerical semigroups starting with $\N$.

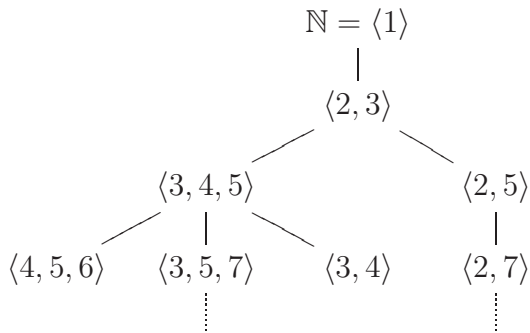
{
    \centerline{
      \xymatrix @R=1pc @C=1pc{
      & & \N=\langle 1\rangle & \\
      & & \langle 2,3\rangle\ar@{-}[u] &\\
      & \langle 3,4,5\rangle \ar@{-}[ur] & & \langle 2,5\rangle \ar@{-}[ul] \\
      \langle 4,5,6\rangle \ar@{-}[ur] & \langle 3,5,7\rangle \ar@{-}[u] \ar@{..}[d] & \langle
      3,4\rangle \ar@{-}[ul]  & \langle 2,7\rangle \ar@{-}[u] \ar@{..}[d]\\
      & & & &      }
      }
}

The level of a vertex in a directed graph is the length of the path connecting this vertex with the root.  Note that in $\Gr(\Sem)$ the level of a vertex coincides with its genus as numerical semigroup.  Therefore, the Bras-Amorós' conjecture quoted in the end of the first section can be reformulated by saying that in   $\Gr(\Sem)$  there are more vertices in the  $(n+1)$th level than in the  $n$th one.

 A Frobenius variety is a nonempty family $\VV$ of numerical semigroups such that
\begin{enumerate}[1)]
\item if $S,T\in \VV$, then $S\cap T\in \VV$,

\item if $S\in \VV$, $S\neq \mathbb N$, then $S\cup \{\F(S)\}\in \VV$.
\end{enumerate}
 The concept of Frobenius variety was introduced in  \cite{variedades-frobenius} with the aim of generalizing most of the results in \cite{patterns, spm, arf, saturated}. In particular, the semigroups that belong to a   Frobenius variety can be arranged as a directed tree with similar properties to those of $\Gr(\Sem)$.

Clearly, $\Sem$ is a Frobenius variety. If $A\subseteq \N$, then $\{S\in \Sem~|~A\subseteq S\}$ is also a  Frobenius variety. In particular, $\Over(S)$, the set of all numerical semigroups that contain  $S$, is a Frobenius variety. We next give some interesting examples of  Frobenius varieties.

Inspired by  \cite{Arf}, Lipman introduces and motivates in \cite{lipman} the study of  Arf rings. The characterization of them via their numerical semigroup of values, brings us to the following concept: a numerical semigroups $S$ is said to be Arf if for every $x,y,z\in S$, with $x,y\geq z$ we have $x+y-z\in S$. It is proved in \cite{arf} that the set of Arf numerical semigroups is a Frobenius variety.

Saturated rings were introduced independently in three distinct ways by  Zariski (\cite{zariski}), Pham-Teissier (\cite{pham}) and Campillo (\cite{campillo}), although the definitions given in these works are equivalent on algebraically closed fields of characteristic zero. Like in the case of numerical semigroups with  the Arf property, saturated numerical semigroups appear when characterizing these rings in terms of their numerical semigroups of values. A numerical semigroup  $S$ is saturated if for every  $s,s_1,\ldots, s_r\in S$ with $s_i\leq s$ for all $i\in \{1,\ldots,r\}$ and $z_1,\ldots,z_r\in \Z$ being integers such that $z_1s_1+\cdots +z_rs_r\geq 0$, then we have $s+z_1s_1+\cdots +z_rs_r\in S$. It is proved in \cite{saturated} that the set of saturated numerical semigroups is a  Frobenius variety.

The class of Arf and Saturated numerical semigroups is also closed under quotients by positive integers as shown in \cite{dobbs-smith}, though the larger class of maximal embedding dimension numerical semigroups is not (if $S$ is a numerical semigroup, then $\e(S)\leq \m(S)$; a numerical semigroup is said to be a maximal embedding dimension semigroup, or to have maximal embedding dimension, if  $\e(S)=\m(S)$). What is the Frobenius variety generated by maximal embedding dimension numerical semigroups?

As a consequence of \cite{spm} and \cite{ns-toms}, it can be deduced that the set of numerical semigroups that admit a Toms decomposition is a  Frobenius variety. Every semigroup with embedding dimension two admits a Toms decomposition. Is the variety of numerical semigroups admitting a Toms decomposition the least  Frobenius variety containing all semigroups with embedding dimension two?

The idea of pattern of a numerical semigroup was introduced in  \cite{patterns} with the aim of trying to generalize the concept of Arf numerical
semigroup. A pattern $P$ of length $n$ is a linear homogeneous polynomial with non-zero integer coefficients in $x_1,\ldots, x_n$ (for $n=0$ the only pattern is $p=0$). We will say that numerical semigroup $S$ admits a pattern  $a_1x_1+\ldots +a_nx_n$ if for every sequence  $s_1\geq s_2\geq \cdots \geq s_n$ of elements in  $S$, we have $a_1s_1+\cdots +a_ns_n\in S$. We denote by $\Sem_P$ the set of all numerical semigroups that admit a pattern  $P$. Then the set of numerical semigroups with  the Arf property is  $\Sem_{x_1+x_2-x_3}$. It is proved in \cite{patterns} that for every pattern $P$ of a special type (strongly admissible), $\Sem_P$ is a Frobenius variety. What varieties arise in this way? It would be interesting to give a weaker definition of pattern such that every variety becomes the variety associated to a pattern.

The intersection of  Frobenius varieties is again a  Frobenius variety. This fact allows us to construct new Frobenius varieties from  known Frobenius varieties and moreover, it allows us to talk of the  Frobenius variety generated by a family  $X$ of numerical semigroups. This variety will be denoted by  $\FV(X)$, and it is defined to be the intersection of all Frobenius varieties containing  $X$. If $X$ is finite, then $\FV(X)$ is finite and it is shown in   \cite{variedades-frobenius} how to compute all the elements of $\FV(X)$.

Let $\VV$ be a Frobenius variety. A submonoid $M$ of $\N$ is a $\VV$-monoid if it can be expressed as an intersection of elements of   $\VV$. It is clear that the intersection of $\VV$-monoids is again a  $\VV$-monoid. Thus given $A\subseteq \N$ we can define the $\VV$-monoid generated by  $A$  as the intersection of all  $\VV$-monoids containing  $A$. We will denote by $\VV(A)$ this $\VV$-monoid and we will say that  $A$ is a $\VV$-system of generators of it. If there is no proper subset of $A$ being a $\VV$-system of generators  $\VV(A)$, then $A$ is a minimal $\VV$-system of generators of  $\VV(A)$. It is proved in \cite{variedades-frobenius} that every  $\VV$-monoid admits an unique  minimal $\VV$-system of generators, and that moreover this system is finite.

We define the directed graph $\Gr(\VV)$ in the same way we defined  $\Gr(\Sem)$, that is, as the graph whose vertices are the elements of  $\VV$, and $(T,S)\in\VV\times \VV$ is an edge of the above graph  if $S=T\cup\{\F(T)\}$. This graph is a tree with root  $\N$ (\cite{variedades-frobenius}). Moreover, the sons of a semigroup  $S$ in $\VV$ are $S\setminus\{x_1\},\ldots, S\setminus\{x_r\}$, where $x_1,\ldots,x_r$ are the minimal $\VV$-generators of $S$ greater than $\F(S)$. This fact allows us to find all the elements of the variety  $\VV$ from $\N$.

 The following figure represents part of the tree associated to the variety of numerical semigroups with the Arf property.

{\tiny
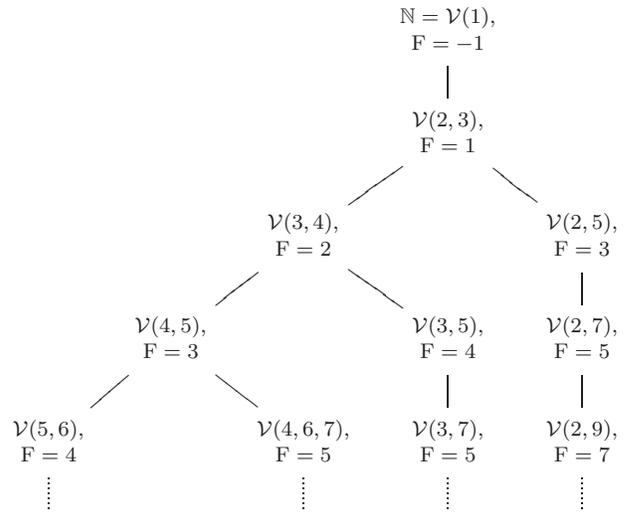
\centerline {\xymatrix @R=1pc @C=1pc{
  &                &                & {\begin{matrix}\N=\VV(1),\\{\rm F}=-1\end{matrix}}\ar@{-}[d] &   \\
  &                &                & {\begin{matrix}\VV(2,3),\\{\rm F}=1\end{matrix}}\ar@{-}[dl]\ar@{-}[dr]&                      \\
  & & {\begin{matrix}\VV(3,4),\\{\rm F}=2\end{matrix}}\ar@{-}[dl]\ar@{-}[dr]
  & & {\begin{matrix}\VV(2,5),\\ {\rm F}=3\end{matrix}}\ar@{-}[d] \\
  & {\begin{matrix}\VV(4,5),\\{\rm F}=3\end{matrix}} \ar@{-}[dl]\ar@{-}[dr] & & {\begin{matrix}\VV(3,5),\\{\rm F}=4\end{matrix}}\ar@{-}[d] & {\begin{matrix}\VV(2,7),\\{\rm F}=5\end{matrix}}\ar@{-}[d] \\
  {\begin{matrix}\VV(5,6),\\{\rm F}=4\end{matrix}}\ar@{-}@{..}[d] & & {\begin{matrix}\VV(4,6,7),\\{\rm F}=5\end{matrix}}\ar@{..}[d] & {\begin{matrix}\VV(3,7),\\{\rm F}=5\end{matrix}}\ar@{..}[d] &  {\begin{matrix}\VV(2,9),\\{\rm F}=7\end{matrix}}\ar@{..}[d] \\
  & & & & }
}}

The following figure represents part of the tree  corresponding to saturated numerical semigroups.

{\tiny 
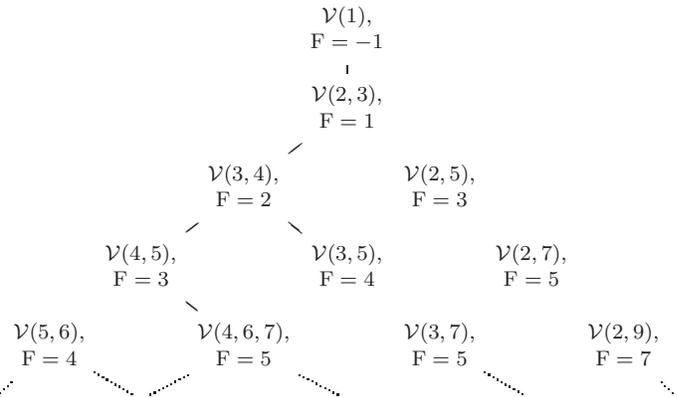
\centerline{\xymatrix @R=.25pc @C=.05pc{
 & & & & {\begin{matrix} \VV(1),\\{\rm F}=-1\end{matrix}} \ar@{-}[d] &  \\
 & & & & {\begin{matrix} \VV(2,3),\\{\rm F}=1\end{matrix}}\ar@{-}[dl]\ar@{-}[dr] & \\
 & & & {\begin{matrix} \VV(3,4),\\{\rm F}=2\end{matrix}}\ar@{-}[dl]\ar@{-}[dr] & &
 {\begin{matrix} \VV(2,5),\\{\rm F}=3\end{matrix}}\ar@{-}[dr] \\
 & & {\begin{matrix} \VV(4,5),\\{\rm F}=3\end{matrix}}\ar@{-}[dl] \ar@{-}[dr] & &
 {\begin{matrix} \VV(3,5),\\{\rm F}=4\end{matrix}} \ar@{-}[dr] & & {\begin{matrix}
 \VV(2,7),\\{\rm F}=5\end{matrix}}\ar@{-}[dr]\\
 & {\begin{matrix} \VV(5,6),\\ {\rm F}=4\end{matrix}} \ar@{..}[dl] \ar@{..}[dr] & &
 {\begin{matrix} \VV(4,6,7),\\ {\rm F}=5\end{matrix}} \ar@{..}[dl] \ar@{..}[dr] & &
 {\begin{matrix} \VV(3,7),\\ {\rm F}=5\end{matrix}} \ar@{..}[dr] & &  {\begin{matrix}
 \VV(2,9),\\ {\rm F}=7\end{matrix}} \ar@{..}[dr]  \\ & & & & & & & & &
 %\\  \vdots & & \vdots ~ \vdots &  & \vdots & & \vdots  & & \vdots
} } }

As a generalization of  Bras-Amorós' conjecture, we can raise the following question. If $\VV$ is a Frobenius variety, does there exist on $\Gr(\VV)$ more vertices in the $(n+1)$th level than  in the   $n$th one? The answer to this question is no, as it is proved in \cite[Example 26]{variedades-frobenius}. However, the same question  in the case of  $\VV$ being infinite remains open. Another interesting question would be characterizing those Frobenius varieties that verify the Bras-Amorós' conjecture.

If $\VV$ is a Frobenius variety and $S\in \VV$, then it is known that $S$ admits an unique minimal $\VV$-system of generators, and moreover it is finite. The cardinality of the set above is an invariant of  $S$ that will be called the embedding $\VV$-dimension of $S$, and it will be denoted by $\e_\VV(S)$. As a generalization of Wilf's conjecture, we would like to characterize those Frobenius varieties  $\VV$ such that for every $S\in \VV$, then $\e_\VV(S)\g(S)\leq (\e_\VV(S)-1)(\F(S)+1)$.

Clearly, the  Frobenius variety generated by irreducible numerical semigroups is  $\Sem$, the set of all numerical semigroups. What is the Frobenius variety generated only by the symmetric ones? and by the  pseudo-symmetric ones?

\section{Presentations of a numerical semigroup}

 Let $(S,+)$ be a commutative monoid. A congruence $\sigma$ over $S$ is an equivalence relation that is compatible with addition, that is,  if $a\sigma b$ with $a,b\in S$, then $(a+c)\sigma (b+c)$ for all $c\in S$. The set $S/\sigma$ endowed with the operation $[a]+[b]=[a+b]$ is a monoid. We will call it the quotient monoid of  $S$ by $\sigma$.

If $S$ is generated by $\{s_1,\ldots,s_n\}$, then the map $\varphi: \N^n\to S,\ (a_1,\ldots, a_n)\mapsto a_1s_1+\cdots +a_ns_n$ is a monoid epimorphism. Therefore $S$ is isomorphic to ${\N^n}/\sim_S$, where $\sim_S$ is the kernel congruence of $\varphi$, that is, $a\sim_S b$ if $\varphi(a)=\varphi(b)$.

The intersection of congruences over a monoid $S$ is again a congruence over  $S$. This fact allows us, given $\sigma\subseteq S\times S$, to define the concept of congruence generated by  $\sigma$ as the intersection of all congruences over  $S$ containing $\sigma$, and it will be denoted by $\langle \sigma\rangle$.

Rédei proves in \cite{redei} that every congruence over $\N^n$ is finitely  generated, that is, there exists a subset of $\N^n\times \N^n $ with finitely many elements generating it. As a consequence we have that giving a finitely generated monoid is, up to isomorphism, equivalent to giving a finite subset of  $\N^n\times \N^n$.

If $S$ is a numerical semigroup with  minimal generators system  $\{n_1,\ldots, n_e\}$, then there exists a finite subset $\sigma$ of $\N^e\times \N^e$ such that $S$ is isomorphic to $\N^e/{\langle \sigma\rangle}$. We say that $\sigma$ is a presentation of $S$. If moreover $\sigma$ has the least possible cardinality, then $\sigma$ is a minimal presentation  of $S$.

A (non directed) graph  $G$ is a pair $(V,E)$, where $V$ is a nonempty set of elements called vertices, and $E$ is a subset of $\{\{u,v\}~|~ u,v\in V, u\neq v\}$. The non ordered pair $\{u,v\}$ will be denoted by $\overline{u v}$, and if it belongs to $E$, then we say that it is an edge of $G$. A sequence of the form $\overline{v_0v_1}, \overline{v_1v_2},\ldots, \overline{v_{m-1}v_m}$ is a path of length $m$ connecting the vertices $v_0$ and $v_m$. A graph is connected if any two distinct vertices are connected by a path. A graph $G'=(V',E')$ is said to be a subgraph  of $G$ if $V'\subseteq V$ and $E'\subseteq E$. A connected component of $G$ is a maximal connected subgraph of $G$. It is well known (see for instance \cite{narsingh}) that a connected graph with $n$ vertices has at least $n-1$ edges. A (finite) tree with  $n$ vertices is a connected graph with $n-1$ edges.

Let us remind now the  method described in \cite{presentaciones-num} for computing the minimal presentation of a numerical semigroup.  Let $S$ be a numerical semigroup with minimal system of generators $\{n_1,\ldots,n_e\}$. For each $n\in S$, let us define $G_n=(V_n,E_n)$, where $V_n=\{n_i~|~ n-n_i\in S\}$ and $E_n=\{\overline{n_in_j}~|~ n-(n_i+n_j)\in S,i\neq j\}$. If $G_n$ is connected, we take $\sigma_n=\emptyset$. If $G_n$ is not connected and  $V_1,\ldots, V_r$ are the sets of vertices corresponding to the  connected components in  $G_n$, then we define $\sigma_n=\{(\alpha_1,\alpha_2),(\alpha_1,\alpha_3),\ldots, (\alpha_1,\alpha_r)\}$, where $\alpha_i\in \varphi^{-1}(n)$ and its $j$-th component is zero whenever  $n_j\not \in V_i$. It is proved in \cite{presentaciones-num} that $\sigma=\bigcup_{n\in S}\sigma_n$ is a minimal presentation for $S$. Let us notice that the set $\textrm{Betti}(S)=\{n\in S~|~ G_n \hbox{ is not connected}\}$ is finite, and that its cardinality is an invariant of $S$. A line of research could be the study of  $\textrm{Betti}(S)$, and its relation with other invariants of $S$ mentioned above. In \cite{betti-unico} affine semigroups (and thus numerical semigroups) with a single Betti element are studied. What are those numerical semigroups having two or three Betti elements?

%particular we could start characterizing those semigroups  $S$ such that $\textrm{A}(S)$ has one element (or two, three,...).
%\textcolor{red}{Para uno já foi feito}

It is also shown in \cite{presentaciones-num} how all the minimal presentations of a semigroup are. In particular, we can determine whether a numerical semigroup admits a unique minimal presentation. Motivated by the idea of generic ideal, we may ask what are the numerical semigroups which admit a unique minimal presentation, and characterize them in terms of their minimal generators.

If $S$ is a numerical semigroup, then the cardinality of a minimal presentation of $S$ is greater than or equal to $\e(S)-1$. Those semigroups that attain this bound are said to be complete intersections. This kind of semigroup has been well studied, and Delorme gives in \cite{delorme} a good characterization of them. Every numerical semigroup with embedding dimension two is a complete intersection, and every complete intersection is symmetric (see \cite{herzog}). We raise the following questions. What semigroups can be expressed as the quotient of a complete intersection by a positive integer? What is the least Frobenius variety containing all the complete intersection numerical semigroups?

Let $S_1$ and $S_2$ be two numerical semigroups minimally generated by  $\{n_1,\ldots, n_r\}$ and $\{n_{r+1},\ldots,n_e\}$, respectively. Let $\lambda \in S_1\setminus \{n_1,\ldots,n_r\}$ and $\mu\in S_2\setminus\{n_{r+1},\ldots, n_e\}$, such that $\gcd\{\lambda,\mu\}=1$. We  then say that $S=\langle \mu n_1,\ldots, \mu n_r, \lambda n_{r+1},\ldots, \lambda n_e\rangle$ is a gluing to  $S_1$ and $S_2$. It is proved in \cite{libro-sn} how given minimal presentations of $S_1$ and $S_2$, one easily gets a minimal presentation of $S$. The characterization given by Delorme in \cite{delorme}, with this notation, can be reformulated in the following way: a numerical semigroup is a complete intersection if and only if is a gluing to two numerical semigroups that are a complete intersection. A consequence of this result is that the set of semigroups that are a complete intersection is the least family of numerical semigroups containing  $\N$ being closed under gluing.  It is well known that the family of numerical symmetric semigroups is also closed under gluing (\cite{libro-sn}). It would be interesting to study other families closed under gluing. Which is the least family containing those semigroups with maximal  embedding dimension and closed under gluing?

Bresinsky gives in  \cite{bresinsky} a family of numerical semigroups with embedding dimension four and with cardinality of its minimal presentations arbitrarily large. This fact proves that the cardinality of  a minimal presentation of a numerical semigroup cannot be  upper bounded  just in function of its embedding dimension.  Bresinski also proves in \cite{bresinsky2} that the cardinality for a minimal presentation of a symmetric numerical semigroup with embedding dimension four can only be three or five. It is conjectured in \cite{rosales-sim-arb} that if $S$ is a numerical semigroup with  $\e(S)\geq 3$, then the cardinality of a minimal presentation for  $S$ is less than or equal to $\frac{\e(S)(\e(S)-1)}2-1$. Barucci  \cite{barucci-na}  proves with the semigroup $\langle 19,23,29,31,37\rangle$ that the conjecture above is not true. However, the problem of determining if the cardinality
of a minimal presentation of a symmetric numerical semigroup  can be bounded in function of the embedding dimension remains open.

Let $\sigma$ be a finite subset of $\N^n\times \N^n$. By using the results in \cite{libro, libro-sn} it is possible to determine  algorithmically whether $\frac{\N^n}{\langle \sigma\rangle}$ is isomorphic to a numerical semigroup. However we miss in the literature families of subsets  $\sigma$ of $\N^n$ so that we can assert, without using algorithms, that $\N^n/{\langle \sigma\rangle}$ is isomorphic to a numerical semigroup. More specifically, we suggest the following problem: given  
\begin{multline*}
\sigma=\{((c_1,0,\ldots,0),(0,a_{1_1},\ldots, a_{1_n})),\ldots, \\
((0,\ldots,c_n),(a_{n_1},\ldots,a_{n_{n-1}},0))\},
\end{multline*} 
which conditions the integers $c_i$ and $a_{j_k}$  have to verify so that $\N^n/{\langle \sigma\rangle}$ is isomorphic to a numerical semigroup? Herzog proved in \cite{herzog} that embedding dimension three numerical semigroups always have a minimal presentation of this form. Neat numerical semigroups introduced by Komeda in \cite{komeda-neat} are also of this form.

%\textcolor{blue}{Komeda para 4 generadores da aquellos que llama neat}

%\textcolor{red}{Falar dos sim\'etricos e pseudo-sim\'etricos de edim 4; tipo?}

\section{Numerical semigroups with embedding dimension three}

Herzog proves in \cite{herzog} that a numerical semigroup with embedding dimension three is symmetric if and only if it is a complete
 intersection. This fact allows us to characterize symmetric numerical semigroups with embedding dimension three in the following way (see   \cite{libro-sn}). A numerical semigroup $S$ with $\e(S)=3$ is symmetric if and only if $S=\langle a m_1, a m_2, b m_1+c m_2\rangle$, with $a$, $b$, $c$, $m_1$ and $m_2$ nonnegative integers, such that
  $m_1$, $m_2$, $a $ and $b+c$ are greater than or equal to two and $\gcd\{m_1,m_2\}=\gcd\{a,b m_1+c m_2\}=1$. Moreover, as it is proved in \cite{libro-sn}, $\F(\langle a m_1, a m_2, b m_1+c m_2\rangle)= a(m_1 m_2-m_1-m_2) + (a-1)(b m_1+c m_2)$. We also have a formula for the genus, since $S$ is symmetric, $\g(S)=\frac{\F(S)+1}2$. Finally, we also know the type, since it is proved in \cite{fgh} that a numerical semigroup is symmetric if and only if its type is equal to one.

We study in \cite{pseudo-sim} the set of pseudo-symmetric numerical semigroups with embedding dimension three.  In particular, we give the following
characterization. A numerical semigroup $S$ with $\e(S)=3$ is pseudo-symmetric if and only if for some ordering of its minimal generators, by taking $\Delta=\sqrt{(\sum n_i)^2-4(n_1n_2+n_1n_3+n_2n_3-n_1n_2n_3)}$, then   $\left\{\frac{n_1-n_2+n_3+\Delta}{2n_1}, \frac{n_1+n_2-n_3+\Delta}{2n_2}, \frac{-n_1+n_2+n_3+\Delta}{2n_3}\right\}\subset \N$. Moreover, in this case, $\F(\langle n_1, n_2, n_3\rangle)=\Delta-(n_1+n_2+n_3)$. We also know the genus and the type, since if $S$ is a pseudo-symmetric numerical semigroups, then $\g(S)=\frac{\F(S)+2}2$ and by \cite{fgh}, $\t(S)=2$.

Bresinsky (\cite{bresinsky}) and Komeda (\cite{komeda-neat}) fully characterize those symmetric and pseudo-symmetric numerical semigroups, respectively, with embedding dimension four. They show that their minimal presentations always have cardinality five.

Curtis proves in \cite{curtis} the impossibility of giving an algebraic formula for the Frobenius number of a numerical semigroup in terms of its minimal generators on embedding dimension three. We raise the following question. Given a polynomial $f(x_1,x_2,x_3,x_4)\in \Q[x_1,x_2,x_3,x_4]$, study the family of numerical semigroups $S$ such that if $S$ is minimally generated by  $n_1< n_2 < n_3$, and $F$ is the Frobenius number of $S$, then $f(n_1,n_2,n_3,F)=0$.

 Our aim now is studying the set of numerical semigroups with embedding dimension three in general.  By \cite{fgh}, we know that these semigroups have type one or two, and by using \cite{johnson, rodseth} if we are concerned with   the Frobenius number and the genus, we can focus ourselves in those numerical semigroups whose minimal generators are pairwise relatively prime. The following result appears in  \cite{tres}. Let $n_1$, $n_2$ and $n_3$ three pairwise relatively prime positive integers. Then the  system of equations
  \[
  \begin{array}{l}
    n_1= r_{12}r_{13}+ r_{12}r_{23} +r_{13}r_{32},\\
    n_2= r_{13}r_{21}+ r_{21}r_{23} +r_{23}r_{31},\\
    n_3= r_{12}r_{31}+ r_{21}r_{32} +r_{31}r_{32}.
  \end{array}
  \]
has a (unique) positive integer solution  if and only if $\{n_1,n_2,n_3\}$ generates minimality $\langle n_1, n_2, n_3\rangle$. In \cite{tres} the authors give formulas for the pseudo-Frobenius number and the genus of $\langle n_1, n_2, n_3\rangle$ from the solutions of the above system. Thus it seems natural to ask, given  positive integers $r_{ij}$, with $i,j\in\{1,2,3\}$, when $r_{12}r_{13}+ r_{12}r_{23} +r_{13}r_{32}$, $r_{13}r_{21}+ r_{21}r_{23}+ r_{23}r_{31}$ and $r_{12}r_{31}+ r_{21}r_{32} +r_{31}r_{32}$ are pairwise relatively prime?

Let $S$ be a numerical semigroup minimally generated by three positive integers  $n_1$, $n_2$ and $n_3$ being pairwise relatively prime. For each $i\in \{1,2,3\}$, let $c_i=\min \{x\in \N\setminus\{0\}~|~ xn_i\in \langle \{n_1,n_2,n_3\}\setminus \{n_i\}\rangle\}$. In \cite{tres} formulas for  $\F(S)$ and $\g(S)$ from  $n_i$ and $c_i$ ($i\in\{1,2,3\}$) are given. Therefore, if we had a formula for computing  $c_3$ from $n_1$ and $n_2$, we would have solved the problems raised by Frobenius for embedding dimension three.  Note that $c_3$ is nothing but the multiplicity of the proportionally modular semigroup  $\frac{\langle n_1,n_2\rangle}{n_3}$. It is proved in \cite{pequeno} that if $u$ is a positive integer such that  $un_2\equiv 1\pmod {n_1}$, then $\frac{\langle n_1,n_2\rangle}{n_3}= \{ x\in \Z~|~ u n_2 n_3 x \bmod n_1 n_2\leq n_3 x\}$. We suggest in this line the problem of finding a formula that allows us to give the multiplicity of $\S(u n_2 n_3,n_1 n_2, n_3)$ from   $n
 _1$, $n_2$ and $n_3$.

 Fermat's Last Theorem  asserts that for any integer  $n\geq 3$, the Diophantine equation $x^n + y^n= z^n$ does not admit an integer solution such that $x y z \neq 0$. As it is well known, this theorem was proved by Wiles, with the help of Taylor, in 1995 (\cite{wiles1, wiles2}) after 300 years of fruitless attempts. Let us observe that for $n\geq 3$, the Diophantine equation $x^n+y^n=z^n$ has no solution verifying $z y z\neq 0$ with some of the factors equal to  $1$. Therefore in order to solve this equation it can be supposed that  $x$, $y$ and $z$ are integers greater than or equal to two, and pairwise relatively prime. It is proved in  \cite{tesis-juan}, that Fermat's Last Theorem
 is equivalent to the following statement:  if $a$, $b$ and $c$ are integers greater than or equal to two,  pairwise relatively prime, and $n$ is an integer greater than or equal to three, then the proportionally modular numerical semigroup  $\frac{\langle a^n, b^n\rangle}c$ is not minimally generated by  $\{a^{n}, c^{n-1}, b^{n}\}$. It would be interesting to prove this fact without using Fermat's last Theorem.
 
%\textcolor{blue}{Mencionar trabajos de Bresinski y Komeda para dimensión de inmersión 4}

\section{Non-unique factorization invariants}

Let $S$ be a numerical semigroup minimally generated by $\{n_1<\cdots < n_e\}$. Then we already know that $S$ is isomorphic to $\mathbb N^e/\sim_S$, where $\sim_S$ is the kernel congruence of the epimorphism $\varphi:\mathbb N^e\to S$, $(a_1,\ldots,a_e)\mapsto a_1n_1+\ldots + a_e n_e$.

For $s\in S$, the elements in $\mathsf Z(s)=\varphi^{-1}(s)$ are known as factorizations of $s$. Given $(x_1,\ldots,x_e)\in \mathbb Z(s)$, its length is $|x|=x_1+\cdots +x_e$. The set of lengths of $s$ is $\mathsf L(s)=\{ |x| ~|~ x\in \mathsf Z(s)\}$. If $\mathsf L(s)=\{ l_1< l_2 < \cdots < l_t\}$, then the set of differences of lengths of factorizations of $s$ is $\mathsf \Delta(s)=\{l_2-l_1,\ldots, l_t-l_{t-1}\}$. Moreover $\Delta(S)=\bigcup_{s\in S}\Delta(s)$. These sets are known to be eventually periodic (\cite{per}).

The elasticity of $s\in S$ is $\rho(s)=\frac{\max\mathsf L(s)}{\min\mathsf L(s)}$, and $\rho(S)=\sup_{s\in S}(\rho(s))$, which turns out to be a maximum (\cite{atomic}). For numerical semigroups it is well known that $\rho(S)=\frac{n_e}{n_1}$. 

For $x=(x_1,\ldots, x_e),y=(y_1,\ldots, y_e)\in \mathbb N^e$, the greatest common divisor of $x$ and $y$ is $\gcd(x,y)=(\min(x_1,y_1),\ldots, \min(x_e,y_e))$. The distance between $x$ and $y$ is $\mathrm d(x,y)= \max\{|x-\gcd(x,y)|,|y-\gcd(x,y)|\}$. 

An $N$-chain (with $N$ a positive integer) joining two factorizations $x$ and $y$ of $s\in S$ is a sequence $z_1,\ldots, z_t$ of factorizations of $s$ such that $z_1=x$, $z_t=y$ and $\mathrm d(z_i,z_{i+1})\le N$. The catenary degree of $s$, $\mathsf c(s)$, is the least $N$ such that for every two factorizations $x$ and $y$ of $s$, there is an $N$-chain joining them. The catenary degree of $S$ is $\mathsf c(S)=\sup_{s\in S}\{\mathsf c(s)\}$. This supremum is a maximum and actually $\mathsf c(S)=\max_{s\in \mathrm{Betti}(S)}\mathsf c(s)$ (\cite{catenary-fg}). It was asked by F. Halter-Koch whether this invariant is periodic, that is, if there exists $n\in S$ such that for $s$ ``big enough'', $\mathsf c(s+n)=\mathsf c(s)$.

The tame degree of $s\in S$, $\mathsf t(s)$, is the minimum $N$ such that for any $i\in \{1,\ldots,e\}$ with $s-n_i\in S$ and any $x\in \mathsf Z(s)$, there exists $y=(y_1,\ldots, y_e)$, such that $y_i\neq 0$ and $\mathsf d(x,y)\le N$. The tame degree of $S$ is $\mathsf t(S)=\sup_{s\in S}(\mathsf t(S))$. This supremum is again a maximum and it is reached in the (finite) set of elements of the form $n_i+w$ with $w\in S$ such that $w-n_j\not\in S$ for some $j\neq i$. F. Halter-Koch also proposed the problem of studying the eventual periodicity of $S$.

The invariant $\omega(S,s)$ is the least positive integer such that whenever $s$ divides $s_1+\cdots+s_k$ for some $s_1,\ldots, s_k\in S$, then $s$ divides $s_{i_1}+\cdots+s_{i_{\omega(S,s)}}$ for some  $\{i_1,\ldots,i_{\omega(S,s)}\}\subseteq \{1,\ldots, k\}$. The $\omega$-primality of $S$ is defined as $\omega(S)=\max\{\omega(S,n_1),\ldots, \omega(S,n_e)\}$. In \cite{omega} it is highlighted that numerical semigroups fulfilling $\omega(S)\neq \mathsf t(S)$ are rare. A characterization for numerical semigroups fulfilling this condition should be welcomed.

Another problem proposed by A. Geroldinger is to determine when can we find a numerical semigroup and an element in it with a given set of lengths.

%\textcolor{red}{Actualizar bibl}


\begin{thebibliography}{10}

\bibitem{Arf} C. Arf, Une interpr\'etation alg\'ebrique de la suite des ordres
  de multiplicit\'e d'une branche alg\'ebrique, Proc. London Math. Soc.,
  \textbf{20} (1949), 256-287.

\bibitem{barucci-na} V. Barucci, Numerical semigroup algebras, in  Multiplicative ideal
theory in commutative algebra,  39--53, Springer, New York, 2006.

\bibitem{barucci} V. Barucci, D. E. Dobbs, M. Fontana, Maximality
  Properties in Numerical Semigroups and Applications to One-Dimensional
  Analytically Irreducible Local Domains, Memoirs of the Amer.  Math. Soc.
  {\bf 598} (1997).
  
\bibitem{omega} V. Blanco, P. A. García-Sánchez, A. Geroldinger, Semigroup-theoretical characterizations of arithmetical invariants with applications to numerical monoids and Krull monoids, Illinois J. Math, to appear.

\bibitem{b-a-conj} M. Bras-Amorós, Fibonacci-like behavior of the number of numerical semigroups of a given genus, Semigroup Forum, \textbf{76} (2008), 379--384.

\bibitem{patterns} M. Bras-Amorós, P. A. García-Sánchez, Patterns on numerical
semigroups, Linar Alg. Appl. {\bf 414} (2006), 652-669.

\bibitem{bresinsky} H. Bresinsky, On prime ideals with generic zeo $x_i=t^{n_i}$, Proc. Amer. Math. Soc. {\bf 47} (1975), 329-332.

\bibitem{bresinsky2} H. Bresinsky, Symmetric semigroups of integers generated by four elements, Manuscripta Math. \textbf{17} (1975), 205-219.

\bibitem{bullejos} M. Bullejos, J. C. Rosales, Proportionally  modular  diophantine  inequalities  and the Stern-Brocot tree, Mathematics of Computation \textbf{78} (2009), 1211-1226

\bibitem{campillo} A. Campillo, On saturation of curve singularities (any
  characteristic), Proc. of Symp. in Pure Math. \textbf{40} (1983), 211-220.

\bibitem{catenary-fg} S. T. Chapman, P. A. García-Sánchez, D. Llena, V. Ponomarenko, and J. C. Rosales, The catenary and tame degree in finitely generated commutative cancellative monoids, Manuscripta Math. \textbf{120} (2006), 253-264.

\bibitem{per} S. T. Chapman, R. Hoyer and N. Kaplan, Delta Sets of Numerical Monoids are Eventually Periodic,
 Aequationes Math. \textbf{77}(2009), 273--279.

\bibitem{curtis} F. Curtis,  On formulas for the Frobenius number of
a numerical semigroup.  Math. Scand. {\bf 67}  (1990),  no. 2, 190-192.

\bibitem{spm} M. Delgado, P. A. García-Sánchez, J. C. Rosales, J. M.
Urbano-Blanco, Systems of proportionally modular Diophantine
inequalities, Semigroup Forum, \textbf{76} (2008), 469-488.

\bibitem{delorme} C. Delorme, Sous-mono\"{\i}des d'intersection compl\`ete de
  $\N$, Ann. Scient. \'Ecole Norm. Sup. (4), {\bf 9} (1976), 145-154.

\bibitem{dobbs} D. E. Dobbs, G. L. Matthews, On a question of Wilf concerning numerical semigroups, en Focus on Commutative Rings Research, 193-202, Nova Sci. Publ., New York, 2006.

\bibitem{dobbs-smith} D. E. Dobbs, H. J. Smith, Numerical semigroups whose fractions are of maximal embedding dimension, Semigroup Forum \textbf{82} (2011), no. 3, 412-422

\bibitem{fgh} R. Fr\"oberg, C. Gottlieb, R. H\"aggkvist, On numerical
  semigroups, Semigroup Forum, {\bf 35} (1987), 63-83.
  
\bibitem{betti-unico} P. A. Garc\'{\i}a-Sánchez, I. Ojeda, J. C. Rosales, Affine semigroups having a unique Betti element, to appear in Journal of Algebra and its Applications.

\bibitem{gkp} R. L. Graham, D. E. Knuth, O. Patashnik, Concrete mathematics. A foundation for computer science. Addison-Wesley Publishing Company, Advanced Book Program, Reading, MA, 1989.

\bibitem{herzog} J. Herzog, Generators and relations of abelian semigroups and semigroup rings, Manuscripta Math. {\bf 3} (1970), 175-193.

\bibitem{johnson} S. M. Johnson, A linear Diophantine problem, Can. J. Math.,
  {\bf 12} (1960), 390-398.

\bibitem{kaplan} N. Kaplan, Counting numerical semigroups by genus and some
cases of a question of Wilf, Journal of Pure and Applied Algebra, \textbf{216} (2012) 1016-1032.


\bibitem{komeda-neat} J. Komeda, On the existence of Weierstrass points with a certain semigroup generated by 4 elements, Tsukuba J. Math, \textbf{6} (1982)m 237-270.

\bibitem{komeda} J. Komeda, Non-weierstras numerical semigroups, Semigroup Forum \textbf{57}(1998), 157-185.

\bibitem{kunz} E. Kunz, The value-semigroup of a one-dimensional Gorenstein
  ring, Proc. Amer. Math. Soc., {\bf 25} (1973), 748-751.

\bibitem{lipman} J. Lipman, Stable ideals and Arf rings, Amer. J. Math.,
  \textbf{93} (1971), 649-685.

\bibitem{narsingh} D. Narsingh, Graph Theory with Applications to Engineering
and Computer Science, Prentice Hall Series in Automatic Computation, 1974.

\bibitem{pham} F. Pham, B. Teissier, Fractions lipschitziennes et
  saturations de Zariski des alg\'ebres analytiques complexes, Centre Math.
  École Polytech., Paris, 1969.  Actes du Congrès International des
  Mathématiciens (Nice, 1970), Tome 2, pp. 649--654.  Gauthier-Villars, Paris,
  1971.

\bibitem{alfonsin} J. L. Ram\'{\i}rez Alfons\'{\i}n, The Diophantine Forbenius
Problem, Oxford University Press, 2005.

\bibitem{redei} L. R\'{e}dei, The theory of finitely generated commutative
  semigroups, Pergamon, Oxford-Edinburgh-New York, 1965.

\bibitem{aureliano} A. M. Robles-Pérez, J. C. Rosales, Equivalent proportionally modular Diophantine inequalities, Archiv der Mathematik \textbf{90} (2008), 24-30.

\bibitem{dobles} A. M. Robles-Pérez, J. C. Rosales y P. Vasco, The doubles of a numerical semigroup, J. Pure Appl. Algebra (aceptado para publicación).

\bibitem{rodseth} Ö. J. Rödseth, On a linear Diophantine problem of Frobenius, J. Reine Angew. Math. \textbf{301} (1979), 431-440.

\bibitem{presentaciones-num} J. C. Rosales,  An algorithmic method to compute a
minimal relation for any numerical semigroup, Internat. J. Algebra
Comput. {\bf 6} (1996), no. 4, 441-455.

\bibitem{rosales-sim-arb} J. C. Rosales, Symmetric numerical semigroups with arbitrary
multiplicity and embedding dimension, Proc. Amer. Math. Soc. {\bf 129} (8)
(2001), 2197-2203.

\bibitem{quotient-psim} J. C. Rosales, One half of a pseudo-symmetric numerical semigroup, Bull. London Math. Soc. \textbf{40} (2008), 347-352.

\bibitem{variedades-frobenius} J. C. Rosales,
Families of numerical semigroups closed under finite intersections and for the
Frobenius number, Houston J. Math. \textbf{34} (2008), 339-348.


\bibitem{interseccion-sim} J. C. Rosales, M. B. Branco, Numerical semigroups that can be
expressed as an intersection of symmetric numerical semigroups, J. Pure Appl.
Algebra {\bf 171} (2-3) (2002), 303--314.

\bibitem{irreducibles-pacific} J. C. Rosales, M. B. Branco, Irreducible numerical
semigroups, Pacific J. Math. {\bf 209} (2003), 131-143.

%
\bibitem{libro} J. C. Rosales, P. A. Garc\'{\i}a-S\'{a}nchez, Finitely
  generated commutative monoids, Nova Science Publishers, New York, 1999.
%

\bibitem{tres} J. C. Rosales, P. A. García-Sánchez, Numerical semigroups with embedding
dimension three, Archiv der Mathematik {\bf 83} (2004), 488-496.

\bibitem{pseudo-sim} J. C. Rosales, P. A. García-Sánchez, Pseudo-symmetric numerical semigroups
with three generators, J. Algebra {\bf 291} (2005), 46-54.

\bibitem{quotient-sym} J. C. Rosales, P. A. García-Sánchez,
Every numerical semigroup is one half of infinitely
many symmetric numerical semigroups, Comm. Algebra \textbf{36} (2008), 2910-2916.

%
\bibitem{libro-sn} J. C. Rosales, P. A. Garc\'{\i}a-S\'{a}nchez, Numerical
    semigroups, Springer, 2009.
%

\bibitem{ns-toms} J. C. Rosales, P. A. García-Sánchez,
Numerical semigroups having a Toms decompositoin, Canadian Math. Bull. \textbf{51} (2008), 134-139.


\bibitem{atomic} J. C. Rosales, P. A. García-Sánchez, J. I. García-García, Atomic commutative monoids and their elasticity, Semigroup Forum \textbf{68} (2004), 64-86.

\bibitem{arf} J. C. Rosales, P. A. García-Sánchez, J. I. García-García,
M. B. Branco, Arf numerical semigroups, J. Algebra {\bf 276} (2004), 3-12.

\bibitem{saturated} J. C. Rosales, P. A. García-Sánchez, J. I. García-García, M. B.
Branco, Saturated numerical semigroups, Houston J. Math. {\bf 30} (2004),
321-330.

\bibitem{oversemigroups} J. C. Rosales, P. A. García-Sánchez, J. I. García-García, J. A.
Jimenez-Madrid,  The oversemigroups of a numerical semigroup, Semigroup Forum
{\bf 67} (2003), 145-158.

\bibitem{fg} J.C. Rosales, P. A. García-Sánchez, J.I. García-García, J.A. Jiménez-Madrid,
 Fundamental gaps in numerical semigroups, J. Pure Appl. Alg. {\bf 189} (2004),
 301-313.

\bibitem{proportionally} J. C. Rosales, P. A. García-Sánchez, J. I. García-García, J. M.
Urbano-Blanco, Proportionally modular Diophantine inequalities, J. Number
Theory {\bf 103} (2003), 281-294.

\bibitem{modular} J. C. Rosales, P. A. García-Sánchez, J. M. Urbano-Blanco,
Modular Diophantine inequalities and numerical semigroup, Pacific J. Math.
{\bf 218} (2) (2005), 379-398.

\bibitem{tres-anos} J. C. Rosales, P. A. García-Sánchez, J. M. Urbano-Blanco,
The set of solutions of a proportionally modular Diophantine inequality,
J. Number Theory. \textbf{128} (2008), 453-467.

\bibitem{opened} J. C. Rosales, J. M. Urbano-Blanco,
Opened modular numerical semigroups,  J. Algebra  {\bf 306}  (2006), 368-377.

\bibitem{full} J. C. Rosales, J.  M. Urbano-Blanco, Proportionally modular Diophantine inequalities and full semigroups,  Semigroup Forum  {\bf 72}  (2006),  362-374.

\bibitem{contraidas} J. C. Rosales y J. M. Urbano-Blanco, Contracted modular Diophantine inequalities and numerical semigroups, Math.Inqual. Appl. 10 (2007), 491-498.

\bibitem{pequeno} J. C. Rosales, P. Vasco, The smallest positive integer that is solution of a proportionally modular Diophantine inequality, Math. Inequal. Appl. \textbf{11} (2008), 203- 212.

\bibitem{alessio} A. Sammartano, Numerical semigroups with large embedding
  dimension satisfy Wilf's conjecture, Semigroup Forum. To appear.

\bibitem{sylvester} J. J. Sylvester, Mathematical questions with their
  solutions, Educational Times {\bf 41} (1884), 21.
  
\bibitem{swanson} I. Swanson, Every numerical semigroup is one over $d$ of infinitely many
   symmetric numerical semigroups, in Commutative algebra and its applications, Walter de Gruyter, Berlin, 2009, 383--386.
   
\bibitem{toms} A. Toms, Strongly perforated $K\sb 0$-groups of simple $C\sp
*$-algebras, Canad. Math. Bull.  {\bf 46}  (2003), 457-472.

\bibitem{tesis-juan} J. M. Urbano-Blanco, Semigrupos numéricos proporcionalmente modulares, Tesis Doctoral, Universidad de Granada, Spain, March 2005.

\bibitem{wiles1} A. Wiles, Modular elliptic curves and Fermat's last theorem, Ann. Math. \textbf{141} (1995), 443-551.

\bibitem{wiles2} A. Wiles y R. Taylos, Ring-Theoretic properties of certain Hecke algebras, Ann. Math. \textbf{141} (1005), 553-572.

\bibitem{wilf} H. S. Wilf, Circle-of-lights algorithm for money changing
problem, Am. Math. Mo.  {\bf 85} (1978), 562-565.

\bibitem{zariski} O. Zariski, General theory of saturation and saturated local
  rings I, II, III, Amer. J. Math. {\bf 93} (1971), 573-684, 872-964,
  {\bf 97} (1975), 415-502.

\bibitem{zhai} A. Zhai, Fibonacci-like growth of numerical semigroups with a
  given genus. arXiv:1111.3142v1.

\bibitem{zhao} Y. Zhao, Constructing numerical semigroups of a given genus Semigroup Forum, \textbf{80} (2010) 242-254

\end{thebibliography}
\end{document}